\documentclass[11pt]{amsart}

\usepackage{latexsym}
\usepackage{amscd}
\usepackage[all]{xy}
\usepackage[mathscr]{euscript}

\normalsize

\RequirePackage{xspace}
\newcommand{\thought}[1]{}
\renewcommand{\thought}[1]{ \textbf{[#1]}}

\usepackage{enumerate}        % Fancy enumerations;  used for Roman numbering.
\newenvironment{roenumerate}{\begin{enumerate}[\upshape (i)]}{\end{enumerate}}
\newcommand\nc {\newcommand}
\newcommand\rnc{\renewcommand}

\newcount\blopone
\newcount\xone
\newcount\xtwo
\newcount\ytwo

\newtheorem{theorem}{Theorem}[section]
\newtheorem{prop}[theorem]{Proposition}
\newtheorem{com}[theorem]{Comment}
\newtheorem{apl}[theorem]{Application}
\newtheorem{exercise}[theorem]{Exercise}
\newtheorem{redu}[theorem]{Reduction}
\newtheorem{refinement}[theorem]{Refinement}
\newtheorem{summary}[theorem]{Summary}
\newtheorem{importnota}[theorem]{Important Notation}
\newtheorem{prblm}[theorem]{Problem}
\newtheorem{notation}[theorem]{Notation}
\newtheorem{explanation}[theorem]{Explanation}
\newtheorem{defin}[theorem]{Definition}
\newtheorem{caution}[theorem]{Caution}
\newtheorem{remark}[theorem]{Remark}
\newtheorem{reminder}[theorem]{Reminder}
\newtheorem{illustration}[theorem]{Illustration}
\newtheorem{observation}[theorem]{Observation}
\newtheorem{lemma}[theorem]{Lemma}
\newtheorem{construction}[theorem]{Construction}
\newtheorem{discussion}[theorem]{Discussion}
\newtheorem{corollary}[theorem]{Corollary}
\newtheorem{example}[theorem]{Example}
\newtheorem{conclusion}[theorem]{Conclusion}
\newtheorem{sketch}[theorem]{Sketch}
\newtheorem{triviality}[theorem]{Triviality}
\newtheorem{proto}[theorem]{Prototype Quasifibration}
\newtheorem{cauex}[theorem]{Cautionary Example}
\newtheorem{hypo}[theorem]{Hypothesis}
\newtheorem{subth}{ }[theorem]
\newtheorem{case}{Case}[theorem]
\newtheorem{ssubth}{ }[subth]
\newtheorem{facts}[theorem]{Facts}
\newtheorem{history}[theorem]{Historical Survey}
\newtheorem{proofs}[theorem]{Discussion of the Proofs, Old and New}

\nc\tri[1]{\begin{triviality}
\label{#1}}
\nc\fac[1]{\begin{facts}
\label{#1}
\begin{em}}
\nc\app[1]{\begin{apl}
\label{#1}
\begin{em}}
\nc\skt[1]{\begin{sketch}
\label{#1}
\begin{em}}
\nc\hst[1]{\begin{history}
\label{#1}
\begin{em}}
\nc\pfs[1]{\begin{proofs}
\label{#1}
\begin{em}}
\nc\cas[1]{\begin{case}
\label{#1}
\begin{em}}
\nc\rfn[1]{\begin{refinement}
\label{#1}}
\nc\prt[1]{\begin{proto}
\label{#1}}
\nc\lem[1]{\begin{lemma}
\label{#1}}
\nc\pro[1]{\begin{prop}
\label{#1}}
\nc\thm[1]{\begin{theorem}
\label{#1}}
\nc\dis[1]{\begin{discussion}
\label{#1}
\begin{em}}
\nc\cor[1]{\begin{corollary}
\label{#1}}
\nc\dfn[1]{\begin{defin}
\label{#1}}
\nc\sthm[1]{\begin{subth}
\label{#1}}
\nc\exm[1]{\begin{example}
\label{#1}
\begin{em}}
\nc\obs[1]{\begin{observation}
\label{#1}
\begin{em}}
\nc\plm[1]{\begin{prblm}
\label{#1}
\begin{em}}
\nc\rmk[1]{\begin{remark}
\label{#1}
\begin{em}}
\nc\rmd[1]{\begin{reminder}
\label{#1}
\begin{em}}
\nc\ntn[1]{\begin{notation}
\label{#1}
\begin{em}}
\nc\exe[1]{\begin{exercise}
\label{#1}
\begin{em}}
\nc\xpl[1]{\begin{explanation}
\label{#1}
\begin{em}}
\nc\smr[1]{\begin{summary}
\label{#1}
\begin{em}}
\nc\cau[1]{\begin{caution}
\label{#1}
\begin{em}}
\nc\hyp[1]{\begin{hypo}
\label{#1}}
\nc\imn[1]{\begin{importnota}
\label{#1}
\begin{em}}
\nc\rdn[1]{\begin{redu}
\label{#1}
\begin{em}}
\nc\cax[1]{\begin{cauex}
\label{#1}
\begin{em}}
\nc\cmt[1]{\begin{com}
\label{#1}
\begin{em}}
\nc\con[1]{\begin{construction}
\label{#1}
\begin{em}}
\nc\ill[1]{\begin{illustration}
\label{#1}
\begin{em}}
\nc\ssthm[1]{\begin{ssubth}
\label{#1}
\begin{em}}
\nc\cnc[1]{\begin{conclusion}
\label{#1}
\begin{em}}

\nc\elem{\end{lemma}}
\nc\erdn{\end{em}\end{redu}}
\nc\erfn{\end{refinement}}
\nc\eprt{\end{proto}}
\nc\ethm{\end{theorem}}
\nc\ecor{\end{corollary}}
\nc\edfn{\end{defin}}
\nc\esthm{\end{subth}}
\nc\epro{\end{prop}}
\nc\etri{\end{triviality}}
\nc\eexm{\end{em}
\end{example}}
\nc\eobs{\end{em}
\end{observation}}
\nc\ecmt{\end{em}
\end{com}}
\nc\efac{\end{em}
\end{facts}}
\nc\eapp{\end{em}
\end{apl}}
\nc\ermk{\end{em}
\end{remark}}
\nc\ermd{\end{em}
\end{reminder}}
\nc\eill{\end{em}
\end{illustration}}
\nc\eplm{\end{em}
\end{prblm}}
\nc\ecas{\end{em}
\end{case}}
\nc\eskt{\end{em}
\end{sketch}}
\nc\ecau{\end{em}
\end{caution}}
\nc\ecax{\end{em}
\end{cauex}}
\nc\eimn{\end{em}
\end{importnota}}
\nc\entn{\end{em}
\end{notation}}
\nc\eexe{\end{em}
\end{exercise}}
\nc\expl{\end{em}
\end{explanation}}
\nc\edis{\end{em}
\end{discussion}}
\nc\econ{\end{em}
\end{construction}}
\nc\esmr{\end{em}
\end{summary}}
\nc\ehst{\end{em}
\end{history}}
\nc\epfs{\end{em}
\end{proofs}}
\nc\ehyp{
\end{hypo}}
\nc\ecnc{\end{em}
\end{conclusion}}
\nc\essthm{\end{em}
\end{ssubth}}

\nc\sst{\scriptstyle}
\newcommand{\comment}[1]{}
\newcommand{\ri}{\longrightarrow}

\newcommand{\D}{{\mathbf D}}

\nc\op{^{\hbox{\rm\tiny op}}}
\nc\mth{^{\hbox{\rm\tiny th}}}

\nc\script{\mathscr}
\nc\z{\zeta}
\nc\bc{{\mathbb{BC}}}
\nc\ct{{\script T}}
\nc\cf{{\script F}}
\nc\cg{{\script G}}
\nc\ch{{\script H}}
\nc\ck{{\script K}}
\nc\cl{{\script L}}
\nc\cv{{\script V}}
\nc\ce{{\script E}}
\nc\cs{{\script S}}
\nc\car{{\script R}}
\nc\cd{{\script D}}
\nc\cc{{\script C}}
\nc\ca{{\script A}}
\nc\ci{{\script I}}
\nc\cj{{\script J}}
\nc\co{{\script O}}
\nc\cu{{\script U}}
\nc\cx{{\script X}}
\nc\Cp{{\script P}}
\nc\cq{{\script Q}}
\nc\cy{{\script Y}}
\nc\cz{{\script Z}}
\nc\bd{\begin{description}}
\nc\ed{\end{description}}
\nc\ctob{{\script C}at\big(\ci^{op},\ca\big)}
\nc\clim{{\ds\mathop{\rm lim}_{\ds\longleftarrow}}\,}
\nc\climi{\clim_{\!i}\,}
\nc\climn{\clim^{\!n}\,}
\nc\colim{{\ds\mathop{\rm colim}_{\ds\la}}}
\nc\colimj{{\ds\mathop{\rm colim}_{\ds\la}}{}_{j\,}}
\nc\oa{\overline{\ca}}
\nc\s{\sigma}
\nc\ta{\tau}
\nc\os{\overline\sigma}
\nc\ot{\overline\tau}
\nc\T{\Sigma}
\nc\Tm{\Sigma^{-1}}
\nc\de[1]{{\mathop{\rm deg(#1)}}}
\nc\Ad[1]{\mathop{\rm Ad}(#1)}
\nc\ad[1]{\mathop{\rm ad}(#1)}
\nc\kth{{\it K}--theory}
\nc\loc[1]{{\text{\rm Loc(#1)}}}
\nc\coloc[1]{{\text{\rm Coloc}(#1)}}

\def\der #1 {D\left(#1\right)}
\nc\prf{\begin{proof}}
\nc\eprf{\end{proof}}
\nc\ds{\displaystyle}
\nc\Tor{\text{\rm Tor}}

\nc\cb{{\script B}}
\nc\ab{{\script A}b}

\nc\be{\begin{roenumerate}}
\nc\ee{\end{roenumerate}}

\nc\cat[1]{{\script C}at\Big({\big\{#1\big\}}\op\,\,,\,\,\ab\Big)}
\nc\csab{{\script C}at\big(\cs^{op},\ab\big)}
\nc\ctab{{\script C}at\Big({\{\ct^\alpha\}}^{op},\ab\Big)}
\nc\csex{{\script E}x\big(\cs^{op},\ab\big)}
\nc\ctex{{\script E}x\Big({\{\ct^\alpha\}}^{op},\ab\Big)}
\nc\sub{\qquad\subset\qquad}
\nc\ctr[1]{{\left.\ct\left(-,#1\right)\right|}_{\cs}}
\nc\ctrf[2]{{\left.\ct\left(#1,#2\right)\right|}_{\cs}}
\nc\Ctr[1]{{\left.\ct\left(-,#1\right)\right|}_{\ct^\alpha}}
\nc\Ctrf[2]{{\left.\ct\left(#1,#2\right)\right|}_{\ct^\alpha}}

\nc\la{\longrightarrow}
\nc\nin{\noindent}
\nc\cad[1]{\text{card}(#1)}
\nc\eq{\quad=\quad}
\nc\BA{\begin{array}{c}}
\nc\EA{\end{array}}
\nc\barr{
\[
\begin{array}{cccccccccccccccc}
}
\nc\earr{
\end{array}
\]
}
\nc\as[1]{{\langle S\rangle}^{#1}}
\nc\sh{\text{\it shift}}

\nc\yy[1]{{\left.\ct\left(-,#1\right)\right|}_{\ct^c}}
\nc\vrep[2]{{\left.\ct\left(#1,#2\right)\right|}_{\ct^\alpha}}
\nc\da{\downarrow}
\nc\Hom{{\mathop{\rm Hom}}}
\nc\HHom{{\script H}{\mathop{\rm om}}}
\nc\End{{\mathop{\rm End}}}
\nc\Ext{{\mathop{\rm Ext}}}
\nc\mr\Modtc
\nc\PExt{{\mathop{\rm PExt}}}
\nc\stm{\text{\rm stmod}(kG)}
\nc\stM{\text{\rm StMod}(kG)}
\nc\e{\varepsilon}
\nc\p{\varphi}

\nc\rs{\s^{-1}A}
\nc\br{{\{\s^{-1}A\}}}
\nc\ra\ri
\nc\y[1]{\mathbf{y}#1}
\nc\x[1]{\mathbf{z}#1}
\nc\mmod[1]{#1\text{--\rm mod}}
\nc\Mod[1]{#1\text{--\rm Mod}}
\nc\Md {\ensuremath{\mathop{\textup{Mod}}}}
\rnc\mod[1]{\ensuremath{\mathop{#1\textup{--mod}}}\xspace}
\nc\Modtc{\Mod{\ct^c}}
\nc\pgldim[1]{\mathop{\rm pgldim}\,#1}
\nc\tf{{\rm [TR5]}}
\nc\tfs{{\rm [TR5$^*$]}}
\nc\Fun{\text{\rm Funct}(F\op,\ab)}
\nc\sym{\text{\rm Sym}}
\nc\sgn{\text{\rm sgn}}
\nc\Pro{\text{\rm Prod}^{}_\alpha(F\op,\ab)}
\nc\Yt[1]{{\left.\Hom_\ct^{}\left(-,#1\right)\right|}_F^{}}
\nc\dl{\delta}
\nc\Proj[1]{#1\text{--\rm Proj}}
\nc\proj[1]{#1\text{--\rm proj}}
\nc\Flat[1]{#1\text{--\rm Flat}}
\nc\Inj[1]{#1\text{--\rm Inj}}
\nc\Ima{\mathrm{Im}}
\nc\Ker{\mathrm{Ker}}
\nc\ov{\overline}
\nc\wt{\widetilde}
\nc\wh{\widehat}
\nc\ph{\varphi}
\nc\tstr{{\it t}--structure}
\nc\spec[1]{{\text{\rm Spec}(#1)}}

\nc\EProd{\text{\rm EProd}}
\nc\ECoprod{\text{\rm ECoprod}}
\nc\Prod{\text{\rm Prod}}
\nc\ldimp{\text{\rm LDim}^{\prod}}
\nc\ldimc{\text{\rm LDim}^{\coprod}}
\nc\gen[2]{{\langle#1\rangle}^{}_{#2}}
\nc\genu[3]{{\langle#1\rangle}^{[#3]}_{#2}}
\nc\ogen[1]{\ov{\langle#1\rangle}}
\nc\ogenun[2]{\ov{\langle#1\rangle}_{#2}^{}}
\nc\ogenu[3]{\ov{\langle#1\rangle}^{[#3]}_{#2}}
\nc\ogenul[3]{\ov{\langle#1\rangle}^{(-\infty,#3]}_{#2}}
\nc\ogenuf[3]{\ov{\langle#1\rangle}^{[#3,\infty)}_{#2}}
\nc\genuf[3]{{\langle#1\rangle}^{[#3,\infty)}_{#2}}
\nc\genul[3]{{\langle#1\rangle}^{(-\infty,#3]}_{#2}}
\nc\dperf[1]{\D^{\mathrm{perf}}(#1)}
\nc\dcoh{\mathbf{D}^b_{\mathrm{coh}}}

\nc\dmcoh{\mathbf{D}^-_{\mathrm{coh}}}
\nc\dscoh{\mathbf{D}^{}_{\mathrm{coh}}}
\nc\RHHom{{\script{RH}}{\mathrm{om}}}
\nc\Coprod{\mathrm{Coprod}}
\nc\COprod{\mathrm{coprod}}
\nc\add{\mathrm{add}}
\nc\Add{\mathrm{Add}}
\nc\Smr{\mathrm{smd}}
\nc\id{\mathrm{id}}
\nc\LL{\mathbf{L}}
\nc\R{\mathbf{R}}
\nc\wi{\wt{\text{\it\i}}}
\nc\exal{\ce\text{\it x}(\ct^\alpha,\ab)}
\nc\exalz{\ce\text{\it x}_{\aleph_0}^{}(\ct^\alpha,\ab)}

\nc\hoco{
\begin{picture}(40,10)
\put(20,0){\makebox(0,0)[b]{\text{\rm Hocolim}}}
\put(5,-2){\vector(1,0){30}}
\end{picture}\,\,}

\nc\holim{
\begin{picture}(40,10)
\put(20,0){\makebox(0,0)[b]{\text{\rm Holim}}}
\put(35,-2){\vector(-1,0){30}}
\end{picture}}

\begin{document}

\author{Amnon Neeman}\thanks{The research was partly supported 
by the Australian Research Council}
\address{Centre for Mathematics and its Applications \\
        Mathematical Sciences Institute\\
        Building 145\\
        The Australian National University\\
        Canberra, ACT 2601\\
        AUSTRALIA}
\email{Amnon.Neeman@anu.edu.au}

\title{The {\it t--}structures generated by objects}

\begin{abstract}
  Let $\ct$ be a well generated triangulated category, and let $S\subset\ct$ be
  a set of objects. We prove that there is a {\it t--}structure on $\ct$
  with $\ct^{\leq0}=\ogenul S{}0$.

  This article is an improvement on the main result of
  Alonso, Jerem{\'{\i}}as and
  Souto~\cite{Alonso-Jeremias-Souto03}, in which the theorem was proved
  under the assumption that $\ct$ has a nice enough model. It should be
  mentioned that the result in \cite{Alonso-Jeremias-Souto03} has been
  influential---it turns out to be interesting to study all of these
  {\it t--}structures.
\end{abstract}



\maketitle

\tableofcontents

\setcounter{section}{-1}

\section{Introduction}
\label{S0}

The main theorem of Alonso, Jerem{\'{\i}}as and
Souto~\cite{Alonso-Jeremias-Souto03} is about constructing {\it t--}structures.
We remind the reader.

Let $\ct$ be a triangulated category with coproducts, and assume
$\big(\ct^{\leq0},\ct^{\geq0}\big)$ is a \tstr\ on $\ct$. Then it's easy and
classical that $\ct^{\leq0}$ is closed under extensions, direct summands,
coproducts and (positive) suspensions.
One can ask the reverse question. Suppose we manage to somehow
concoct a full subcategory $\cs\subset\ct$, closed under all these
operations---in Keller and
Vossieck~\cite{Keller-Vossieck87,Keller-Vossieck88} and 
Alonso, Jerem{\'{\i}}as and
Souto~\cite{Alonso-Jeremias-Souto03} such an $\cs$ is called a
``cocomplete pre-aisle''.
Given a cocomplete pre-aisle $\cs\subset\ct$, one can wonder if there
is a \tstr\ with $\cs=\ct^{\leq0}$.

Keller and Vossieck~\cite[Section~1]{Keller-Vossieck88} showed that a
(cocomplete) pre-aisle
$\cs$ is equal to $\ct^{\leq0}$ if and only if the inclusion $\cs\la\ct$
has a right adjoint. And Alonso, Jerem{\'{\i}}as and
Souto~\cite{Alonso-Jeremias-Souto03} construct such a right adjoint if
$\cs$ is the
cocomplete pre-aisle ``generated'' by a set of objects in $\ct$,
and if $\ct$ has a sufficiently nice model---a locally presentable,
cofibrantly generated model
category will do.
In this article we prove an improvement, we don't assume $\ct$ has any model.
Our main theorem says

\medskip

\nin
{\bf Theorem~\ref{T2.5}.}\ \ 
\emph{
Let $\ct$ be a well generated triangulated category and let $S\subset\ct$
be a set of objects. Then there is a \tstr\ on $\ct$ with 
$\ct^{\leq0}=\ogenul S{}0$, where $\ogenul S{}0$ is our 
notation for the smallest 
cocomplete pre-aisle in $\ct$ containing $S$.}

\medskip

Until this article the only model-free version of such a theorem assumed
that the objects in $S$ are all compact in $\ct$, see
\cite[Thoerem~A.1]{Alonso-Jeremias-Souto03}.

The proof of the more general
theorem in~\cite{Alonso-Jeremias-Souto03}, where
the objects in $S$ aren't restricted to be compact but the category
$\ct$ is assumed to have a nice model,
hinges on a small-object argument in Quillen's sense.
We should say something about our proof which, needless to say, is totally
different.

Suppose the pre-aisle $\cs$ is an aisle, meaning (in the terminology of
Keller and
Vossieck~\cite{Keller-Vossieck87,Keller-Vossieck88}) that there
is a \tstr\ with $\cs=\ct^{\leq0}$. Then for every object $t\in\ct$ there must
be an object $t^{\leq0}\in\cs$, yielding a homological
functor $H(-)=\Hom\big(-,t^{\leq0}\big)$. The idea of this
paper is to construct the functor $H$ directly,
and then use Brown representability
to exhibit it as $\Hom\big(-,t^{\leq0}\big)$. It turns out that,
as long as we're willing to disregard set-theoretic issues,
the definition
of $H$ is simple enough---the reader can find it at the very beginning
of Section~\ref{S1}. The construction
makes sense in great generality, and it is straightforward to
show that $H:\ct\op\la\mathit{AB}$ is a homological functor respecting products.
Here $\mathit{AB}$ stands for large abelian groups, the collection of elements
might not be a set---it needn't belong to our universe.
The only subtle part is the set-theoretic problem: it is a little
tricky to show that $H(x)$ is a small set for every $x\in\ct$.

Finally we should mention one more result we prove in this article.
But first the historical context: assuming all the objects
in the set $S\subset\ct$ are compact, Keller and
Nicol{\'a}s~\cite[Theorem~A.9]{Keller-Nicolas13} give a refinement of
Alonso Jerem{\'{\i}}as and Souto~\cite[Thoerem~A.1]{Alonso-Jeremias-Souto03}.
Not only do they show that the category $\ogenul S{}0$ is the aisle
of a \tstr, they prove further that every object $x\in\ogenul S{}0$
can be expressed as the homotopy colimit of a countable sequence
$x_1^{}\la x_2^{}\la x_3^{}\la\cdots$, where each
$x_\ell^{}$ is an $\ell$--fold extension
of coproducts of positive suspensions of objects in $S$.
We give an analog
of this which holds
when the objects in $S$ aren't assumed compact---it
requires some notation to state our result precisely,
the reader is referred to
Propsition~\ref{P2.9}.

Assuming the existence of a nice model,
as in the proof of the main theorem of Alonso Jerem{\'{\i}}as and
Souto~\cite{Alonso-Jeremias-Souto03}, does not offer an
alternative approach to
Proposition~\ref{P2.9}. The small object argument of
\cite{Alonso-Jeremias-Souto03} constructs sequences which decidedly
\emph{aren't} countable.

\medskip

\nin
{\bf Acknowledgements.}\ \
The author would like to thank Bernhard Keller, for pointing out
the interesting and relevant results of~\cite{Keller-Nicolas13}, and
Michal Hrbek and Bregje Pauwels for their careful reading of the current
version.

\section{$t$--structures via representability}
\label{S1}

\dfn{D1.1}
Let $\ct$ be a triangulated category,
and let $\cs\subset\ct$ be a full, additive
subcategory. Assume that $\cs$ is a pre-aisle, meaning
\be
\item
$\T\cs\subset\cs$.
\item
  $\cs*\cs=\cs$, meaning if $x\la y\la z$ is a triangle with $x,z\in\cs$ then
  $y\in\cs$.
\item
$\Smr(\cs)=\cs$, in other words $\cs$ contains all direct summands of its
objects.
\setcounter{enumiv}{\value{enumi}}
\ee
For every pair of objects $t,t'\in\ct$ we make the following definitions:
\be
\setcounter{enumi}{\value{enumiv}}
\item
The class $\wt H_\cs^{}(t,t')$ has for its elements the pairs of composable
morphisms
$t\la s\la t'$ with $s\in\cs$.
\setcounter{enumiv}{\value{enumi}}
\ee
And now we define a relation on
$\wt H_\cs^{}(t,t')$.
Given two objects $t,t'\in\ct$, then
$R(t,t')\subset\wt H_\cs^{}(t,t')\times\wt H_\cs^{}(t,t')$ is as follows
\be
\setcounter{enumi}{\value{enumiv}}
\item
Suppose we are given
$(h,h')\in\wt H_\cs^{}(t,t')\times\wt H_\cs^{}(t,t')$,
where
\[
h=\{t\stackrel f\la s\stackrel g\la t'\}\qquad\text{ and }\qquad
h'=\{t\stackrel{f'}\la s'\stackrel{g'}\la t'\}\ .
\]
The pair $(h,h')$ belongs to $R(t,t')$
if 
there is in $\ct$ a commutative diagram
\[\xymatrix@C+30pt@R-5pt{
           & s\ar[dr]\ar@/^1pc/[drr]^-g & & \\
t\ar[ur]^-f\ar[dr]_-{f'} & & s''\ar[r] & t' \\
           & s'\ar[ur]\ar@/_1pc/[urr]_-{g'} & & 
}\]
with $s''\in\cs$.
\setcounter{enumiv}{\value{enumi}}
\ee
\edfn

\obs{O1.2}
In passing we note that, given $(h,h')\in R(t,t')$
with $h=\{t\stackrel f\la s\stackrel g\la t'\}$ and
$h'=\{t\stackrel{f'}\la s'\stackrel{g'}\la t'\}$, then
the diagram of Definition~\ref{D1.1}(v) forces
the equality $gf=g'f'$.
\eobs

\exm{E1.3.5}
Suppose we are given composable morphisms
$t\stackrel e\la s\stackrel f\la s'\stackrel g\la t'$
with $s,s'\in\cs$.  Let
\[
h=\{t\stackrel e\la s\stackrel {gf}\la t'\}\qquad\text{ and }\qquad
h'=\{t\stackrel{fe}\la s'\stackrel{g}\la t'\}\ .
\]
The commutative diagram
\[\xymatrix@C+30pt@R-5pt{
           & s\ar[dr]_-f\ar@/^1pc/[drr]^-{gf} & & \\
t\ar[ur]^-e\ar[dr]_-{fe} & & s'\ar[r]|g & t' \\
           & s'\ar[ur]^-\id\ar@/_1pc/[urr]_-{g} & & 
}\]
shows that $(h,h')\in R(t,t')$.
\eexm

\lem{L1.3}
The $R(t,t')$ of Definition~\ref{D1.1}(v) is an equivalence relation.
\elem

\prf
Since $R(t,t')$ is obviously reflexive and symmetric, what needs
proof is transitivity. Suppose therefore that we are given
three elements $h,h',h''\in\wt H_\cs^{}(t,t')$ with
$(h,h')$ and $(h',h'')$ in $R(t,t')$. If 
\[
h=\{t\stackrel f\la s\stackrel g\la t'\},\qquad
h'=\{t\stackrel{f'}\la s'\stackrel{g'}\la t'\},\qquad
h''=\{t\stackrel{f''}\la s''\stackrel{g''}\la t'\}
\]
then the diagrams exhibiting the fact that
$(h,h')$ and $(h',h'')$ belong to $R(t,t')$
assemble to a
commutative diagram
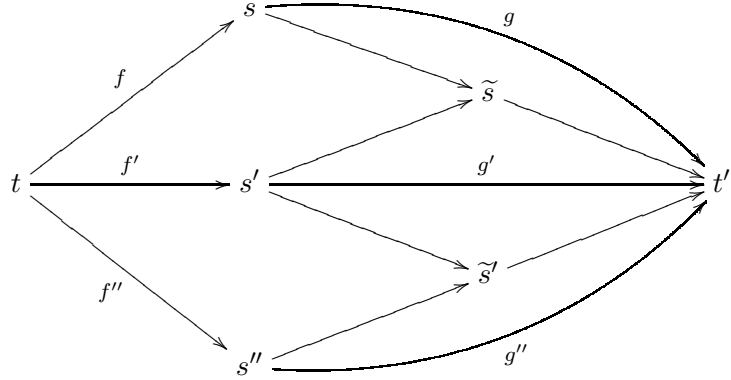
\[\xymatrix @C+50pt@R-5pt{
& s\ar[dr]\ar@/^2pc/[ddrr]^-g & & \\
& & \wt s\ar[rd] &  \\
t\ar[uur]^-f\ar[ddr]_-{f''}\ar[r]^-{f'} & s'\ar[ur]\ar[dr]\ar[rr]^-{g'} & & t' \\
& & \wt s'\ar[ru] &  \\
           & s''\ar[ur]\ar@/_2pc/[uurr]_-{g''} & & 
}\]
The commutative square
\[\xymatrix{
s' \ar[r]\ar[d] & \wt s\ar[d]\\
\wt s'\ar[r] & t'
}\]
may be factored through the homotopy pushout
\[\xymatrix{
s' \ar[r]\ar[d] & \wt s\ar[d]\\
\wt s'\ar[r] & \wt s''\ar[r] & t'
}\]
We remind the reader: by the definition of homotopy pushouts
the square in the last
diagram is commutative, and may be ``folded'' to a triangle
$s'\la \wt s\oplus\wt s'\la\wt s''\la \T s'$.
This makes $\wt s''$ an object in $\cs*\T\cs\subset\cs*\cs=\cs$.
The commutative diagram
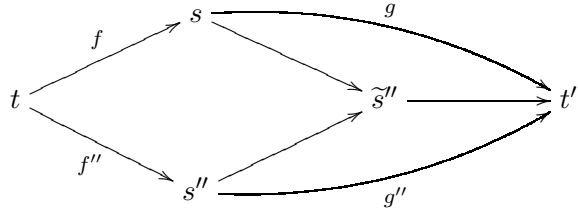
\[\xymatrix@C+30pt@R-5pt{
           & s\ar[dr]\ar@/^1pc/[drr]^-g & & \\
t\ar[ur]^-f\ar[dr]_-{f''} & & \wt s''\ar[r] & t' \\
           & s''\ar[ur]\ar@/_1pc/[urr]_-{g''} & & 
}\]
now establishes that $(h,h'')$ belongs to $R(t,t')$.
\eprf

\dfn{D1.5}
With the notation of Definition~\ref{D1.1}~(iv) and (v), 
we define $H_\cs^{}(t,t')$ to be the quotient of $\wt H_\cs^{}(t,t')$
by the equivalence relation $R(t,t')$.
\edfn

\exm{E1.7}
The example to keep in mind is the following. Suppose $\ct$ is a triangulated
category with a \tstr, and put $\cs=\ct^{\leq0}$. Then an element of
$\wt H_\cs^{}(t,t')$ is a pair of composable morphisms $t\la s\la t'$ with
$s\in\cs=\ct^{\leq0}$, and this composable string may factored further,
uniquely, as
$t\la s\la(t')^{\leq0}\la t'$ where the map $(t')^{\leq0}\la t'$
is the canonical map from the \tstr\ truncation.
Example~\ref{E1.3.5} applies, showing that every element
of $\wt H_\cs^{}(t,t')$
is equivalent to an element $t\stackrel f\la(t')^{\leq0}\la t'$.

Now suppose the elements $t\stackrel f\la(t')^{\leq0}\la t'$ and
$t\stackrel{f'}\la(t')^{\leq0}\la t'$ are equivalent to each other. Then there
must exist in $\ct$ a commutative diagram
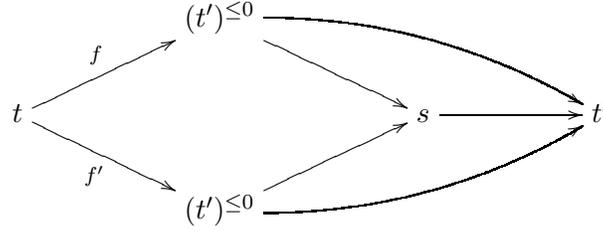
\[\xymatrix@C+30pt@R-5pt{
           & (t')^{\leq0}\ar[dr]\ar@/^1pc/[drr] & & \\
t\ar[ur]^-f\ar[dr]_-{f'} & & s\ar[r] & t' \\
           & (t')^{\leq0}\ar[ur]\ar@/_1pc/[urr] & & 
}\]
with $s\in\cs$. The map $s\la t'$ must also factor as $s\la (t')^{\leq0}\la t'$,
allowing us to replace the above by the commtative diagram
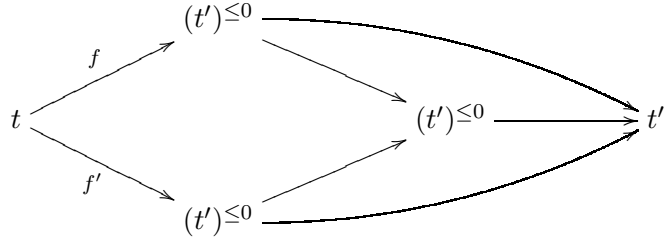
\[\xymatrix@C+30pt@R-5pt{
           & (t')^{\leq0}\ar[dr]\ar@/^1pc/[drr] & & \\
t\ar[ur]^-f\ar[dr]_-{f'} & &(t')^{\leq0} \ar[r] & t' \\
           & (t')^{\leq0}\ar[ur]\ar@/_1pc/[urr] & & 
}\]
But now in the two commutative triangles
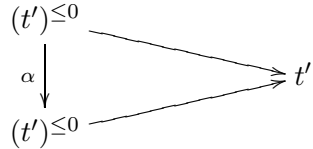
\[\xymatrix@C+50pt@R-20pt{
(t')^{\leq0}\ar[dr] \ar[dd]_\alpha&  \\
   & t' \\
(t')^{\leq0}\ar[ur] & 
}\]
the maps $\alpha$ have to be identities, forcing $f=f'$.

Thus in the special case, where $\cs=\ct^{\leq0}$
for some \tstr, the class
$H_\cs^{}(t,t')$ is canonically isomorphic to $\Hom\big[t,(t')^{\leq0}\big]$.
\eexm

It's natural to wonder how much of the structure, which
$\Hom\big[t,(t')^{\leq0}\big]$
always has, can be constructed on $H_\cs^{}(t,t')$ without knowing
that
$\cs=\ct^{\leq0}$ for some \tstr. This leads us to the next few results.

\dfn{D1.9}
Let $\ct$ be a triangulated category, let $\cs\subset\ct$ be as in
Definition~\ref{D1.1}~(i), (ii) and (iii), and let $H_\cs^{}(t,t')$ be
as
in Definition~\ref{D1.5}. We define the following:
\be
\item
If $h,h'$ are two elements of $H_\cs^{}(t,t')$, then $h+h'$ is
constructed as follows. First choose representatives for 
$h,h'$ in $\wt H_\cs^{}(t,t')$, that is choose
factorizations
\[
t\stackrel{f}\la s \stackrel{g}\la t',\qquad
t\stackrel{f'}\la s' \stackrel{g'}\la t'
\]
which represent the equivalence classes $h,h'$. Then $h+h'$ is the
equivalence
class of 
\[\xymatrix@C+40pt{
\ds t \ar[r]^-{\left(\begin{array}{c}
f \\ f'
\end{array}\right)} & 
s\oplus s' \ar[r]^-{\ds\big(g\,,\,g'\big)} &
t
}\]
\item
The element $0\in H_\cs^{}(t,t')$ is defined to be the equivalence
class of $t\la 0\la t'$.
\ee
\edfn

\nin
The next lemma is straightforward, the proof is left to the reader.

\lem{L1.11}
The operation $+$ of Definition~\ref{D1.9}(i) is well-defined, meaning
the element $h+h'\in  H_\cs^{}(t,t')$ does not depend on the choice of representatives.
The operation $+$ is commutative and associative. The element $0\in
H_\cs^{}(t,t')$ 
of Definition~\ref{D1.9}(ii) satisfies $0+h=h=h+0$. Finally: if $h$ is
represented
by $t\stackrel{f}\la s \stackrel{g}\la t'$ and 
$h'$ is represented by $t\stackrel{f}\la s \stackrel{-g}\la t'$ then $h+h'=0$.

Summarizing: Definition~\ref{D1.9} gives $H_\cs^{}(t,t')$  the
structure of an abelian group.
\elem

\con{C1.13}
Given a morphism $e:x\la y$, we define
$\wt H_\cs^{}(e,t):\wt H_\cs^{}(y,t)\la
\wt H_\cs^{}(x,t)$
to be the map which precomposes with $e$. That is: 
$y\stackrel{f}\la s \stackrel{g}\la t$ goes to $x\stackrel{fe}\la
s \stackrel{g}\la t$.
This map obviously sends equivalent elements of 
$\wt H_\cs^{}(
y,t)$
to equivalent elements of
$\wt H_\cs^{}(x,t)$, hence descends to a map which we will denote
\[\xymatrix@C+50pt{
H_\cs^{}(e,t)\,\,:\,\,H_\cs^{}(y,t)\ar[r] & H_\cs^{}(x,t)
}\]
\econ

\nin
The next lemma is another obvious one.

\lem{L1.15}
Let $t\in\ct$ be an object.
With the notation of Construction~\ref{C1.13}, we have that
$H_\cs^{}(-,t)$ is an additive functor from $\ct$ to abelian
groups. Here
$H_\cs^{}(x,t)$ is an abelian group with the structure given
in Lemma~\ref{L1.11}.
\elem

Slightly more subtle is 

\lem{L1.17}
The functor $H_\cs^{}(-,t)$ of Lemma~\ref{L1.15} is homological.
\elem

\prf
Suppose we are given an object $t\in\ct$, and let
$x\stackrel\alpha\la y\stackrel\beta\la z$ be a
triangle
in $\ct$. We need to show that $H_\cs^{}(-,t)$ takes it to an exact
sequence.

Suppose therefore that we are given an element $h\in H_\cs^{}(y,t)$
which
maps to zero in $H_\cs^{}(x,t)$. Choose a representative for $h$,
meaning
composable morphisms 
$y\stackrel{f}\la
s \stackrel{g}\la t$ with $s\in\cs$. We are given that the image of
$h$ in $H_\cs^{}(x,t)$ vanishes, meaning the composite 
$x\stackrel{f\alpha} \la
s \stackrel{g}\la t$ 
is equivalent to zero. The equivalence means that there  must be a
commutative
diagram
\[\xymatrix@C+30pt@R-5pt{
           & s\ar[dr]_-{g'}\ar@/^1pc/[drr]^-{g} & & \\
x\ar[ur]^-{f\alpha}\ar[dr] & & s'\ar[r]_{g''} & t' \\
           & 0\ar[ur]\ar@/_1pc/[urr] & & 
}\]
with $s'\in\cs$. The composable morphisms
$y\stackrel{f}\la
s \stackrel{g'}\la s' \stackrel{g''}\la t$ and Example~\ref{E1.3.5}
tell us that 
$y\stackrel{f}\la
s \stackrel{g}\la t$,
which is equal to 
$y\stackrel{f}\la
s \stackrel{g''g'}\la t$, is equivalent to
$y\stackrel{g'f}\la
s' \stackrel{g''}\la t$. In other
words $h\in H_\cs^{}(y,t)$ is also represented 
by $y\stackrel{g'f}\la
s' \stackrel{g''}\la t$. But now the vanishing of
the composite
$x\stackrel{\alpha}\la
y\stackrel{g'f}\la
s'$
says that $g'f:y\la s'$ must factor as
$y\stackrel\beta\la z\stackrel\gamma\la s'$. But then
$h\in H_\cs^{}(y,t)$ must be the image of
$z\stackrel\gamma\la s'\stackrel{g''}\la t$ under the map
$H_\cs^{}(\beta,t):H_\cs^{}(z,t)\la H_\cs^{}(y,t)$.
\eprf

\nin
And finally we look at the case where $\ct$ has coproducts.

\lem{L1.19}
Let $\ct$ be a triangulated category with coproducts, let $\cs\subset\ct$ be a
full subcategory closed in $\ct$ under coproducts, and assume that
the hypotheses of Definition~\ref{D1.1}~(i), (ii) and (iii) are satisfied.

Then the functor $H_\cs^{}(-,t)$ of Lemma~\ref{L1.15} respects products.
To expand: we view $H_\cs^{}(-,t)$ as a functor $\ct\op\la\ab$, and
the products in $\ct\op$ are the coproducts in $\ct$. Suppose
$\{x_\lambda^{},\,\lambda\in\Lambda\}$ is a set of objects in $\ct$,
there is always a natural map
\[\xymatrix@C+30pt@R-5pt{
\ds H_\cs^{}\left(\coprod_{\lambda\in\Lambda}x_\lambda^{}\,\,,\,\,t\right)
\ar[r]^-\rho &
\ds \prod_{\lambda\in\Lambda}H_\cs^{}(x_\lambda^{},t)
}\]
and we assert that this map is an isomorphism.
\elem

\prf
The inverse of the canonical map $\rho$ is simple enough.  Given an
element $\prod_{\lambda\in\Lambda}h_\lambda^{}\in\prod_{\lambda\in\Lambda}H_\cs^{}(x_\lambda^{},t)$ choose representatives,
meaning for each $\lambda$ choose  for $h_\lambda^{}$
a representative $x_\lambda^{}\la s_\lambda^{}\la t$. Form
the composite
\[\xymatrix@C+30pt@R-5pt{
\ds \coprod_{\lambda\in\Lambda}x_\lambda^{}
\ar[r]&
\ds \coprod_{\lambda\in\Lambda}s_\lambda^{}
\ar[r] &
t\ .
}\]
The reader can check that this construction gives a well-defined map
\[\xymatrix@C+30pt@R-5pt{
\ds \prod_{\lambda\in\Lambda}H_\cs^{}(x_\lambda^{},t)
\ar[r]^-\s &
\ds H_\cs^{}\left(\coprod_{\lambda\in\Lambda}x_\lambda^{}\,\,,\,\,t\right)\ ,
}\]
meaning the resulting element of
$H_\cs^{}\left(\coprod_{\lambda\in\Lambda}x_\lambda^{}\,\,,\,\,t\right)$
is independent of the choice of representatives.

And now it is an exercise
to check that $\s\rho$ and $\rho\s$ are both identities.
\eprf

\smr{S1.21}
We have learned that, for any $\cs\subset\ct$ as in Definition~\ref{D1.1}
and any object $t\in\ct$,
the functor $H_\cs^{}(-,t)$ of Lemma~\ref{L1.15} is homological. If $\ct$ has
coproducts and $\cs$ is closed in $\ct$ under coproducts then
$H_\cs^{}(-,t)$ also respects products.

Example~\ref{E1.7} showed us that, in the special case where $\cs=\ct^{\leq0}$
for some \tstr, the functor  $H_\cs^{}(-,t)$ is representable---it is naturally
isomorphic to $\Hom\big(-,t^{\leq0}\big)$. The previous paragraph amounts to
saying that, for a general $\cs$, the functor $H_\cs^{}(-,t)$ satisfies
the obvious necessary conditions for representability. 
\esmr

Now we come to

\pro{P1.23}
Suppose $\cs\subset\ct$ are as in Definition~\ref{D1.1}, and let
$t\in\ct$ be an object. Assume the functor $H_\cs^{}(-,t)$
of Lemma~\ref{L1.15} is representable, more explicitly assume
we are given an isomorphism $\Hom(-,x)\la H_\cs^{}(-,t)$ which we fix.
Then
\be
\item
The object $x$ belongs to $\cs\subset\ct$.
\item
There is a unique morphism
$\e:x\la t$, so that the image of $\id:x\la x$ under the map
$\Hom(x,x)\la H_\cs^{}(x,t)$ is represented by
$x\stackrel\id\la x\stackrel\e\la t$.
\item
Any map $s\la t$, with $s\in\cs$, factors in $\ct$ uniquely as
$s\stackrel{\exists!}\la x\stackrel\e\la t$.
\setcounter{enumiv}{\value{enumi}}
\ee
Finally: the subcategory $\cs$ is equal to $\ct^{\leq0}$,
for some \tstr\ on $\ct$,
if and only if $ H_\cs^{}(-,t)$ is representable for every $t\in\ct$.
\epro

\prf
We begin by proving (i), starting with the given isomorphism
$\Hom(-,x)\la H_\cs^{}(-,t)$. The identity map $\id\in\Hom(x,x)$ must go under
the isomorphism $\Hom(x,x)\la H_\cs^{}(x,t)$ to an
element of $H_\cs^{}(x,t)$; choose a representative
$x\stackrel g\la s\stackrel h\la t$ with $s\in\cs$. The composable morphisms
$s\stackrel\id\la s\stackrel h\la t$ may be viewed as
representing an element of $H_\cs^{}(s,t)$,
which must be the image of
an $f\in\Hom(s,x)$ under the isomorphism
$\Hom(s,x)\la H_\cs^{}(s,t)$. We know the image
$x\stackrel g\la s\stackrel h\la t$ of
$\id\in\Hom(x,x)$ under
the isomorphism $\Hom(x,x)\la H_\cs^{}(x,t)$, and Yoneda tells us that
$f$ has to go to $s\stackrel{gf}\la s\stackrel h\la t$. But
$f$ was chosen to have image $s\stackrel\id\la s\stackrel h\la t$, and
therefore
\[
s\stackrel{gf}\la s\stackrel h\la t\qquad\text{ and }\qquad
s\stackrel{\id}\la s\stackrel h\la t
\]
must be equivalent. Precomposing with $g:x\la s$ we have that 
\[
x\stackrel{gfg}\la s\stackrel h\la t\qquad\text{ and }\qquad
x\stackrel{g}\la s\stackrel h\la t
\]
are also equivalent. The natural isomorphism
$\Hom(x,x)\la H_\cs^{}(x,t)$ therefore takes
the elements $fg,\id\in\Hom(x,x)$ to the same image, and we deduce that
$fg=\id$.
Thus
$x$ is a direct summand of $s\in\cs$, and Definition~\ref{D1.1}(iii)
guarantees that $x$ belongs to $\cs$. This proves (i).

To prove (ii) recall that the image of $\id:x\la x$, under
the isomorphism $\Hom(x,x)\la H_\cs^{}(x,t)$, is represented by
the composable pair $x\stackrel f\la s\stackrel g\la t$ which we
chose at the beginning of the proof. But now we know that $x\in\cs$ and
Example~\ref{E1.3.5}, applied to the
composable morphisms
$x\stackrel \id\la x\stackrel f\la s\stackrel g\la t$, tells us that
$x\stackrel f\la s\stackrel g\la t$
is equivalent to  $x\stackrel\id\la x\stackrel{gf}\la t$.
Furthermore this representative is unique: if
$x\stackrel\id\la x\stackrel{h}\la t$ is equivalent to
$x\stackrel\id\la x\stackrel{h'}\la t$ then Observation~\ref{O1.2}
teaches us that $h=h'$.
As in the statement of the current Proposition
this preferred representative, of the image of
$\id\in\Hom(x,x)$ under the isomorphism
$\Hom(x,x)\la H_\cs^{}(x,t)$,
will be
denoted $x\stackrel\id\la x\stackrel{\e}\la t$.

Now for the proof of (iii). Suppose we are given a map $s\stackrel f\la t$
with $s\in\cs$. Then $s\stackrel\id\la s\stackrel{f}\la t$,
viewed as representing an element in $H_\cs^{}(s,t)$,
corresponds under the isomorphism
$H_\cs^{}(s,t)\cong\Hom(s,x)$ to a morphism
$\alpha:s\la x$. The isomorphism
$\Hom(s,x)\la H_\cs^{}(s,t)$ takes $\alpha$ to the equivalence class
of $s\stackrel\alpha\la  x\stackrel{\e}\la t$. The fact
that $s\stackrel\id\la s\stackrel{f}\la t$ is equivalent
to $s\stackrel\alpha\la  x\stackrel{\e}\la t$, coupled
with Observation~\ref{O1.2}, tells us that 
$f=\e\alpha$.

The uniqueness of $\alpha$ is proved as
follows. Suppose we have elements $\alpha,\alpha'$ in $\Hom(s,x)$ with 
$\e\alpha=\e\alpha'$. Consider
the triangle
$s\stackrel{\alpha-\alpha'}\la x\stackrel\beta\la s'\la\T s$. 
As $x,\T s$ both lie in $\cs$ so does $s'$, while
$\e:x\la t$ has to factor  $x\stackrel\beta\la s'\stackrel\gamma\la t$.
The diagram
\[\xymatrix@C+30pt@R-5pt{
           & x\ar[dr]_-{\beta}\ar@/^1pc/[drr]^-\e & & \\
s\ar[ur]^-\alpha\ar[dr]_-{\alpha'} & & s'\ar[r]|\gamma & t \\
           & x\ar[ur]^-{\beta}\ar@/_1pc/[urr]_-{\e} & & 
}\]
now exhibits that the pairs
\[
h=\{s\stackrel\alpha\la x\stackrel\e\la t\},\qquad
h'=\{s\stackrel{\alpha'}\la x\stackrel\e\la t\}
\]
are equivalent as in
Definition~\ref{D1.1}(v). But this means that 
the isomorphism $\Hom(s,x)\cong H_\cs^{}(s,t)$ takes
$\alpha,\alpha'\in\Hom(s,x)$ to the same image, hence $\alpha=\alpha'$.

Now for the last part of the Proposition. Example~\ref{E1.7} teaches us that,
if $\cs=\ct^{\leq0}$ for some \tstr, then $H_\cs^{}(-,t)$ is representable
for every $t\in\ct$. What needs proof is the converse. Assume therefore
that every $H_\cs^{}(-,t)$ is representable. By (i), (ii) and (iii) we
deduce that, for every $t\in\ct$, there exists a morphism $\e:x\la t$ such that
\be
\setcounter{enumi}{\value{enumiv}}
\item
$x\in\cs$
\item
Every morphism $s\la t$, with $s\in\cs$, factors uniquely through $\e$.
\ee
But this exactly says that the inclusion $\cs\la\ct$ has a right adjoint,
and $\e:x\la t$ is the counit of adjunction. From
\cite[Section~1]{Keller-Vossieck88} it now follows that $\cs=\ct^{\leq0}$ for
a (unique) \tstr\ on $\ct$.
\eprf

\dis{D1.25}
Let $\ct$ be a triangulated category with coproducts and
let $S\subset\ct$ be a set of objects. Form the category $\cs=\ogenul S{}0$,
the smallest full subcategory of $\ct$
satisfying the hypotheses of Definition~\ref{D1.1}~(i), (ii) and (iii)
and closed in
$\ct$ under coproducts. Summary~\ref{S1.21} allows us to deduce that
the functor $H_\cs^{}(-,t)$ is homological and respects products---still
in the gorgeous generality of any triangulated category $\ct$ with coproducts.

Now assume $\ct$ is well generated. Then Brown representability holds,
see \cite[Theorem~8.3.3]{Neeman99}---by this theorem showing
that $H_\cs^{}(-,t)$ is representable only requires
solving the set-theoretic problem, we need to
prove that $H_\cs^{}(x,t)$ is a small set for every $x,t\in\ct$.

The proof of Theorem~\ref{T2.5} will show how to solve this set theoretic
problem. And then Proposition~\ref{P1.23} will come to our aid---since
we will know that $H_\cs^{}(-,t)$ is representable for every $t\in\ct$,
Proposition~\ref{P1.23} will allow us to deduce that $\ogenul S{}0=\ct^{\leq0}$
for some \tstr.
\edis

\section{Main theorem}
\label{S2}

\con{C2.1}
Let $\ct$ be a well-generated triangulated category, and let $S\subset\ct$ be
a small set of objects. Choose a
large enough regular cardinal $\alpha$ so that
\be
\item
The category $\ct$ is $\alpha$--compactly generated.
\item  
The set $S$ is contained in $\ct^\alpha$.
\setcounter{enumiv}{\value{enumi}}
\ee
Now we proceed by transfinite induction, on the ordinal
$i\leq\alpha$, to build up full subcategories
$S(i)\subset\ct^\alpha\cap\ogenul S{}0$ as follows:
\be
\setcounter{enumi}{\value{enumiv}}
\item
The objects of $S(1)$ are the coproducts of $<\alpha$ objects in
$\cup_{j=0}^\infty \T^jS$.
\item  
If $i$ is any ordinal $<\alpha$ and $i+1$ is its successor, then
the objects of 
$S(i+1)$ are all the 
coproducts of $<\alpha$ objects in $S(i)*(S(i)$.
\item
If $i'\leq\alpha$ is  limit ordinal, then $S(i')=\cup_{i<i'}S(i)$.
\setcounter{enumiv}{\value{enumi}}
\ee
By induction we see that each $S(i)$ satisfies $\T S(i)\subset S(i)$,
and for any limit ordinal $i$ we have that $S(i)*S(i)=S(i)$.
By the above and the fact that $\alpha$ is a regular cardinal
we have that $S(\alpha)$ satisfies
\be
\setcounter{enumi}{\value{enumiv}}
\item
$S\subset S(\alpha)\subset\ct^\alpha\cap\ogenul S{}0$.
\item
$\T S(\alpha)\subset S(\alpha)$.
\item  
$S(\alpha)*S(\alpha)=S(\alpha)$.
\item
Any 
coproduct of $<\alpha$ objects in $S(\alpha)$ lies in $S(\alpha)$.
\setcounter{enumiv}{\value{enumi}}
\ee
\econ

Now we assert:

\lem{L2.3}
Let the notation be as in Construction~\ref{C2.1}.
Any morphism $t\la s'$, with $t\in\ct^\alpha$ and $s'\in\ogenul S{}0$, can be
factored as $t\la s\la s'$ with $s\in S(\alpha)$.
\elem

\prf
Let $\car\subset\ct$ be the full subcategory defined by the formula
\[
\text{Ob}(\car)\eq
\left\{
r\in\ct\left|\begin{array}{c}
\text{Every morphism $t\la r$, with $t\in\ct^\alpha$,}\\
  \text{can be factored as $t\la s\la r$ with $s\in S(\alpha)$}
\end{array}\right.\right\}  
\]
and we observe the following
\be
\item
$S\subset \car$.
\setcounter{enumiv}{\value{enumi}}
\ee
This is obvious: any map 
$t\la s$ can be factored as
$t\la s\stackrel\id\la s$. 
\be
\setcounter{enumi}{\value{enumiv}}
\item
$\T\car\subset\car$.
\setcounter{enumiv}{\value{enumi}}
\ee
Suppose we are given a morphism $f:t\la\T r$ with $t\in\ct^\alpha$
and $r\in\car$. Then $\Tm f:\Tm t\la r$ is a morphism from $\Tm t\in\ct^\alpha$
to $r\in\car$, and may be factored as 
$\Tm t\la s\la r$ with $s\in S(\alpha)$. Thus $f$ has a factorization as
$t\la\T s\la\T r$ with $\T s\in\T S(\alpha)\subset S(\alpha)$.
\be
\setcounter{enumi}{\value{enumiv}}
\item
The subcategory $\car$ is closed in $\ct$ under coproducts.
\setcounter{enumiv}{\value{enumi}}
\ee
To see this let $\{r_\lambda^{},\,\lambda\in\Lambda\}$ be any
set of objects in $\car$, let $t$ be an object in
$\ct^\alpha$, and suppose we are given a morphism
\[\xymatrix@C+40pt{
t\ar[r]^-f & \ds\coprod_{\lambda\in\Lambda} r_\lambda^{}\ .
}\]
The fact that $t$ belongs to $\ct^\alpha$ permits us to find a
subset $\Lambda'\subset\Lambda$ of cardinality $<\alpha$, and
a factorization
\[\xymatrix@C+40pt{
t\ar[r] & \ds\coprod_{\lambda\in\Lambda'} t_\lambda^{}
\ar[r]^-{\coprod_{\lambda\in\Lambda'} f_\lambda^{}} &
\ds\coprod_{\lambda\in\Lambda'} r_\lambda^{}
\ar@{^{(}->}[r] &  \ds\coprod_{\lambda\in\Lambda} r_\lambda^{}
}\]
with each $t_\lambda^{}\in\ct^\alpha$. But then for each
$\lambda\in\Lambda'$ we can factor $f_\lambda^{}:t_\lambda^{}\la r_\lambda^{}$
as $t_\lambda^{}\la s_\lambda^{}\la r_\lambda^{}$ with $s_\lambda^{}\in S(\alpha)$,
giving a factorization of $f$ as
\[\xymatrix@C+40pt{
t\ar[r] & \ds\coprod_{\lambda\in\Lambda'} s_\lambda^{}
\ar[r] &  \ds\coprod_{\lambda\in\Lambda} r_\lambda^{}\ .
}\]
But since the cardinality of $\Lambda'\subset\Lambda$ is $<\alpha$
and each $s_\lambda^{}$ belongs to $S(\alpha)$,
Construction~\ref{C2.1}(ix) gives that $\coprod_{\lambda\in\Lambda'} s_\lambda^{}$
must belong to $S(\alpha)$.
\be
\setcounter{enumi}{\value{enumiv}}
\item
The subcategory $\car$ satisfies $\car*\car\subset\car$.
\setcounter{enumiv}{\value{enumi}}
\ee
To prove (iv) apply \cite[Lemma~1.5]{Neeman17},
with $\ca=\cc=S(\alpha)$ and $\cx=\cz=\car$; we have that any pair
of morphisms $t\la x$ and $t\la z$, with $t\in\ct^\alpha$ and with $x\in\cx$
and $z\in\cz$, factor (respectively) as $t\la a\la x$ and $t\la c\la z$ with
$a\in\ca$ and $c\in\cc$. Since $S(\alpha)=\cc\subset\ct^\alpha$ and
$\ct^\alpha$ is triangulated \cite[Remark~1.6]{Neeman17} applies.
We conclude that any map $t\la y$, with $t\in\ct^\alpha$ and
$y\in\cx*\cz=\car*\car$, must factor as $t\la b\la y$ with
$b\in\ca*\cc=S(\alpha)*S(\alpha)\subset S(\alpha)$.

Since $\ogenul S{}0$ is the smallest subcategory of $\ct$ satisfying (i),
(ii), (iii) and (iv) it must be contained in $\car$, proving the lemma.
\eprf

\thm{T2.5}
Let $\ct$ be a well generated triangulated category, and let $S\subset\ct$
be a set of objects. Then there is a \tstr\ on $\ct$ with
$\ct^{\leq0}=\ogenul S{}0$.
\ethm

\prf
Put $\cs=\ogenul S{}0$. By Discussion~\ref{D1.25}
it suffices to show that, for every pair of
objects $x,t\in\ct$, the collection $H_\cs^{}(x,t)$ is a set.

For this choose a regular cardinal $\alpha$ large enough so that $\ct$
is $\alpha$--compactly generated, and $S\cup\{x\}\subset\ct^\alpha$. Let
$S(\alpha)$ be as in Construction~\ref{C2.1}. Suppose we are given a
representative $x\stackrel f\la s\stackrel g\la t$ for an element
$h\in H_\cs^{}(x,t)$. Then
$f:x\la s$ is a morphism from $x\in\ct^\alpha$ to $s\in\cs=\ogenul S{}0$,
and Lemma~\ref{L2.3} allows us to factor $f$ as
$x\stackrel{f'}\la s'\stackrel{f''}\la s$
with $s'\in S(\alpha)$. The string of composable morphisms
$x\stackrel{f'}\la s'\stackrel{f''}\la s\stackrel g\la t$,
coupled with Example~\ref{E1.3.5},
yields an equivalence between
$x\stackrel{f''f'}\la s\stackrel{g}\la t$ and
$x\stackrel{f'}\la s'\stackrel{gf''}\la t$, that is
$x\stackrel f\la s\stackrel g\la t$ is equivalent to
$x\stackrel{f'}\la s'\stackrel{gf''}\la t$. Thus every $h\in H_\cs^{}(x,t)$
may be represented as $x\la s'\la t$ with $s'\in S(\alpha)$, and there
is only a set of these.

Note that the equivalence can also be checked without going to very
large objects. If  $x\stackrel{f}\la s\stackrel{g}\la t$ and
$x\stackrel{f'}\la s'\stackrel{g'}\la t$ are equivalent,
and $s$ and $s'$ both lie in $S(\alpha)$, then
the definition of equivalence says there exists a
commutative diagram
\[\xymatrix@C+30pt@R-5pt{
           & s\ar[dr]\ar@/^1pc/[drr]^-g & & \\
x\ar[ur]^-f\ar[dr]_-{f'} & & \wt s\ar[r] & t \\
           & s'\ar[ur]\ar@/_1pc/[urr]_-{g'} & & 
}\]
with $\wt s\in\cs=\ogenul S{}0$. The commutative square
\[\xymatrix{
x\ar[r]^-f\ar[d]_-{f'} & s\ar[d] \\
s'\ar[r] & \wt s
}\]
may be factored through the homotopy pushout, obtaining
\[\xymatrix{
x\ar[r]^-f\ar[d]_-{f'} & s\ar[d] \\
s'\ar[r] & y \ar[r] & \wt s
}\]
The fact that we have a homotopy pushout means there is a triangle
$x\la s\oplus s'\la y\la\T x$, and as $x$ and $s\oplus s'$ belong to $\ct^\alpha$
so does $y$. The map $y\la\wt s$ is a morphism from $y\in\ct^\alpha$
to $\wt s\in\ogenul S{}0$, and Lemma~\ref{L2.3} permits us
to factorize it as $y\la s''\la\wt s$ with $s''\in S(\alpha)$. We obtain
a commutative diagram
\[\xymatrix@C+30pt@R-5pt{
           & s\ar[dr]\ar@/^1pc/[drrr]^-g & & \\
x\ar[ur]^-f\ar[dr]_-{f'} & & s''\ar[r] & \wt s\ar[r] & t \\
           & s'\ar[ur]\ar@/_1pc/[urrr]_-{g'} & & 
}\]
where $s,s',s''$ all lie in $S(\alpha)$, and the objects $x$ and $t$ are given
and fixed. If we delete $\wt s$ the diagram still exhibits the equivalence
of $x\stackrel{f}\la s\stackrel{g}\la t$ and
$x \stackrel{f'}\la s'\stackrel{g'}\la t$, but now there is only a set of
possible such diagrams to consider.

Summarizing: we have proved that $H_\cs^{}(x,t)$ is a set, and the theorem
follows.
\eprf

\rmd{R2,7}
Before the next result we should remind the reader of  basic
notation.
Let $\ct$ be a triangulated category with coproducts and let
$T\subset\ct$ be a set of objects. 
Then, inductively on the integer $\ell>0$,
we construct in $\ct$
full subcategories
$\Coprod_\ell^{}(T)$
as follows:
\be
\item
For $\ell=1$ we declare the  category $\Coprod_1^{}(T)$
to have for objects all possible coproducts of the objects
belonging to $T$.
\item
Assume $\Coprod_\ell^{}(T)$ has been constructed.
Then the category $\Coprod_{\ell+1}^{}(T)$ is defined by
the formula
\[
\Coprod_{\ell+1}^{}(T)\eq
\Coprod_{\ell}^{}(T)*
\Coprod_{1}^{}(T)\ .
\]
This formula means that an object $y\in\ct$ belongs to
$\Coprod_{\ell+1}^{}(T)$ if there exists in $\ct$ a
triangle
$x\la y\la z$ with
$x\in \Coprod_{\ell}^{}(T)$ and
$z\in \Coprod_{1}^{}(T)$.
\ee
\ermd

And now we are ready for the next result, which is an analog
of Keller and Nicol{\'a}s~\cite[Theorem~A.9]{Keller-Nicolas13}.

\pro{P2.9}
Let $\ct$ be a well-generated triangulated category and assume
$S\subset\ct$
is a set of objects. Pick a regular cardinal $\alpha$ so that
\be
\item
$\ct$ is $\alpha$--compactly generated.
\item
$S\subset\ct^\alpha$.
\setcounter{enumiv}{\value{enumi}}
\ee
And now let $S(\alpha)$ be as in Construction~\ref{C2.1}.

Given any object $x\in\ogenul S{}0$, there exists in $\ct$ a 
countable sequence
$x_1^{}\la x_2^{}\la x_3^{}\la\cdots$ such that
\be
\setcounter{enumi}{\value{enumiv}}
\item
$x_\ell^{}$ belongs to $\Coprod_{\ell}^{}\big[S(\alpha)\big]$.
\item
$x\cong\hoco x_\ell^{}$.
\setcounter{enumiv}{\value{enumi}}
\ee
\epro

\prf
Fix for the proof an object $x\in\ogenul S{}0$.
We begin by constructing the sequence $x_1^{}\la x_2^{}\la x_3^{}\la\cdots$ and
a map $\ph:\hoco x_\ell^{}\la x$, and then we will prove
that $\ph$ is an isomorphism.
The construction of the sequence is as follows.
\be
\setcounter{enumi}{\value{enumiv}}
\item
Let $\Lambda_1$ be the set of all maps
$f_\lambda^{}:s_\lambda^{}\la x$, with $s_\lambda^{}\in S(\alpha)$. Define 
\[
x_1^{}\eq\coprod_{\lambda\in\Lambda_1}s_\lambda^{}\ .
\]
The formula makes it clear that 
$x_1^{}$ belongs to $\Coprod_1^{}\big[S(\alpha)\big]$.
We let the map $\ph_1^{}:x_1^{}\la x$ be the obvious.
\item
Assume $x_\ell^{}\in \Coprod_\ell^{}\big[S(\alpha)\big]$ 
and the map $\ph_\ell^{}:x_\ell^{}\la x$ have been defined. Complete
$\ph_\ell^{}$ to a
triangle $y\la x_\ell^{}\stackrel{\ph_\ell^{}}\la x$,
and let $\Lambda_\ell$ be the set
of all maps $\Tm s_\lambda^{}\la y$ with $s_\lambda^{}\in S(\alpha)$.
The composite
$\coprod_{\Lambda_\ell}\Tm s_\lambda^{}\la x_\ell^{}
\stackrel{\ph_\ell^{}}\la x$
factors through the composite $y\la x_\ell^{}\stackrel{\ph_\ell^{}}\la x$
and therefore vanishes, allowing us to construct the commutative
diagram
\[\xymatrix@C+30pt{
\ds\coprod_{\Lambda_\ell}\Tm s_\lambda^{}\ar[r] & x_\ell^{}
\ar[r]^-{g_\ell^{}}\ar[dr]_-{\ph_\ell^{}} &
x_{\ell+1}^{}\ar[r]\ar[d]^-{\ph_{\ell+1}^{}} & \ds\coprod_{\Lambda_\ell}s_\lambda^{}\\
 & & x & \\
}\]
where the row is a triangle. In other words we factor
$\ph_\ell^{}:x_\ell^{}\la x$
as a composite $x_\ell^{}\stackrel{g_\ell^{}}\la x_{\ell+1}^{}
\stackrel{\ph_{\ell+1}^{}}\la x$, and the triangle in the top row
exhibits $x_{\ell+1}^{}$ as belonging to
$\Coprod_\ell^{}\big[S(\alpha)\big]*\Coprod_1^{}\big[S(\alpha)\big]=\Coprod_{\ell+1}^{}\big[S(\alpha)\big]$.
\setcounter{enumiv}{\value{enumi}}
\ee
We have now defined the sequence
$x_1^{}\stackrel{g_1^{}}\la x_2^{}\stackrel{g_2^{}}\la x_3^{}\stackrel{g_3^{}}\la \cdots$
as well as, for each $\ell>0$, a map $\ph_\ell^{}:x_\ell^{}\la x$.
Since the $\ph_\ell^{}$'s are compatible with the $g_\ell^{}$'s,
we may factor the map from the sequence through some $\ph:\hoco x_\ell^{}\la x$.
It remains to prove that any such $\ph$ is an isomorphism.
To this end we prove
\be
\setcounter{enumi}{\value{enumiv}}
\item
Let $t$ be any object in $\ct^\alpha$. 
Then the functor $\Hom(t,-)$ takes the map $\ph_1^{}:x_1^{}\la x$ to an
epimorphism.
\setcounter{enumiv}{\value{enumi}}
\ee
To prove (vii) let $f:t\la x$ be any morphism. Then $f$ is a map from
$t\in\ct^\alpha$ to $x\in\ogenul S{}0$, and Lemma~\ref{L2.3} permits us
to factor $f$ as $t\la s\la x$ with $s\in S(\alpha)$. 
But then $f$ factors through the coproduct of all $s\la x$, which
in (v) was defined to be $\ph_1^{}:x_1^{}\la x$.

Next we assert
\be
\setcounter{enumi}{\value{enumiv}}
\item
Let $t$ be any object in $\ct^\alpha$. Applying the functor 
$\Hom(t,-)$ to the displayed diagram 
in (vi) we obtain
\[
\xymatrix@C+50pt{
\Hom(t,x_\ell^{})\ar[r]^-{\Hom(t,g_\ell^{})} \ar[dr]_{\Hom(t,\ph_\ell^{})}&
\Hom(t,x_{\ell+1}^{})\ar[d]^{\Hom(t,\ph_{\ell+1}^{})}\\
 &\Hom(t,x)
}\]
and we assert that the horizontal map and the slanted map
have the same kernel.
\setcounter{enumiv}{\value{enumi}}
\ee
To prove the assertion note first that the kernel of the
horizontal map is obviously contained in the kernel of the slanted
one, what needs proof is the reverse inclusion.

Choose therefore any map $f:t\la x_\ell^{}$ such that the composite
$t\stackrel f\la x_\ell^{} \stackrel{\ph_\ell^{}}\la x$ vanishes.
Recalling the triangle $y\la x_\ell^{} \stackrel{\ph_\ell^{}}\la x$
of (vi), the map $f$ must factor as $t\la y\la x_\ell^{}$.
But the triangle $\Tm x\la y\la x_\ell^{}$ exhibits $y$ as an
object of $\Tm\ogenul S{}0*\ogenul S{}0\subset\Tm\ogenul S{}0$,
making the morphism $t\la y$ a map from $t\in\ct^\alpha$ to
$y\in\Tm\ogenul S{}0$. Lemma~\ref{L2.3}, applied to the suspension
of this map, allows us to factor $t\la y$  as a composite
$t\la \Tm s\la y$ with
$s\in S(\alpha)$. Hence $t\la y$ factors through the coproduct
of all $\Tm s\la y$, and the triangle in the top row
of the display diagram in (vi) shows that the composite
$t\stackrel f\la x_\ell^{} \stackrel{g_\ell^{}}\la x_{\ell+1}^{}$ must
vanish. We have proved (viii).

It remains to deduce the Proposition from (vii) and (viii), but this is
standard. We remind the reader.

Let $\exal$ be the abelian category, whose objects
are the additive functors
$(\ct^\alpha)\op\la\ab$ which respect products of $<\alpha$
objects. This category is developed extensively in
\cite{Neeman99}, see also the clever alternative description in
Krause~\cite{Krause01}. Let $\cy:\ct\la\exal$ be the functor
taking the object $t\in\ct$
to the object $\cy(t)=\Hom_\ct^{}(-,t)\big|_{\ct^\alpha}^{}$
in the category $\exal$. It's obvious
that $\cy$ is a homological functor, and a theorem that
it respects coproducts---see \cite[Proposition~6.2.6]{Neeman99}.

Now apply the functor $\cy$ to everything in sight. By (vii) the map
$\cy(\ph_1^{}):\cy(x_1^{})\la\cy(x)$ is an epimorphism and,
since this epimorphism factors as
$\cy)(x_1^{})\la\cy(x_\ell^{})\stackrel{\cy(\ph_\ell^{})}\la\cy(x)$
for every integer $\ell>0$,
the maps $\cy(\ph_\ell^{}):\cy(x_\ell^{})\la\cy(x)$
must all be epimorphisms. Now by
(viii) the map $\cy(g_\ell^{}):\cy(x_\ell^{})\la\cy(x_{\ell+1}^{})$
factors canonically as
$\cy(x_\ell^{})\stackrel{\cy(\ph_\ell^{})}
\la\cy(x)\la \cy(x_{\ell+1}^{})$. Therefore the sequence
\[\xymatrix@C+20pt{
\cy(x_1^{}) \ar[r]^-{\cy(g_1^{})} &
\cy(x_2^{}) \ar[r]^-{\cy(g_2^{})} &
\cy(x_3^{}) \ar[r]^-{\cy(g_3^{})} &
\cdots
}\]
is ind-isomorphic to the sequence 
\[\xymatrix@C+20pt{
\cy(x) \ar[r]^-{\id} &
\cy(x) \ar[r]^-{\id} &
\cy(x) \ar[r]^-{\id} &
\cdots
}\]
allowing us to compute, as (for example) in the proof of
\cite[Theorem~8.3.3]{Neeman99}, that the map
$\cy(\ph):\cy\big(\hoco x_\ell^{}\big)\la\cy(x)$ must be an
isomorphism.

But $\alpha$ was chosen large enough so that $\ct^\alpha$ generates,
hence the map $\ph:\hoco x_\ell^{}\la x$ must also be an isomorphism.
\eprf

\def\cprime{$'$}
\providecommand{\bysame}{\leavevmode\hbox to3em{\hrulefill}\thinspace}
\providecommand{\MR}{\relax\ifhmode\unskip\space\fi MR }
\providecommand{\MRhref}[2]{%
  \href{http://www.ams.org/mathscinet-getitem?mr=#1}{#2}
}
\providecommand{\href}[2]{#2}

\end{document}